\documentclass[twoside,leqno]{article}
\usepackage{amsfonts}
\usepackage{amsmath}
\usepackage{amssymb}
\usepackage{amsthm}
\usepackage{graphicx}
\usepackage{epsfig}
\usepackage{amsmath,amssymb,amsfonts}
\usepackage{lineno}

\def\bjga#1{\def\thefootnote{}\footnote{
    \hspace*{-.54cm} $~~~~~$The first author is the corresponding author.}
    \addtocounter{footnote}{-1}\def\thefootnote{\arabic{footnote}}}
\def\mymaketitle#1{\date{\bjga{#1}}\begin{document}\maketitle\thispagestyle{empty}}

\def\babs{\begin{abstractv}}
\def\eabs{\end{abstractv}}

\theoremstyle{definition}
\newtheorem{definition}{Definition}[section]
\newtheorem{example}[definition]{Example}
\newtheorem{remark}[definition]{Remark}
\newtheorem{problem}[definition]{Problem}

\theoremstyle{theorem}
\newtheorem{theorem}{Theorem}[section]
\newtheorem{proposition}[theorem]{Proposition}
\newtheorem{corollary}[theorem]{Corollary}
\newtheorem{lemma}[theorem]{Lemma}

\def\pas{\par\smallskip}
\def\pasn{\par\smallskip\noindent}
\def\pam{\par\medskip}
\def\pamn{\par\medskip\noindent}
\def\pab{\par\bigskip}

\def\header#1#2#3#4#5{
    \markboth{#3}{#4}
    \title{#5}
    \author{#3}
    \date{}
    \mymaketitle{#1-#2}}
\newenvironment{abstractv}{\begin{quote}{\bf Abstract.\ }}{\end{quote}}
\def\msc{{\bf M.S.C. 2010}:\ }
\def\kwd{\\{\bf Key words}:\ }
\def\aua{\par\noindent{\em Author's address:}\pam\noindent}
\def\auas{\par\noindent{\em Authors' addresses:}\pam\noindent}
\def\auac{\par\noindent{\em Authors' address:}\pam\noindent}
\def\bece{\begin{center}}
\def\eece{\end{center}}
\def\bebi{\begin{thebibliography}{40}\setlength{\parskip}{-1mm}}
\def\eebi{\end{thebibliography}}
\def\bibi#1{\bibitem{#1}}

\newcommand{\R}{\mathbb{R}}
\newcommand{\Rt}{\mbox{{\em $\R$}}}
\newcommand{\Rs}{\mbox{\tiny{\R}}}
\newcommand{\C}{\mathbb{C}}
\newcommand{\Cs}{\mbox{\tiny{\C}}}
\newcommand{\Q}{\mathbb{Q}}
\newcommand{\Z}{\mathbb{Z}}
\newcommand{\Zs}{\mbox{\tiny{\Z}}}
\newcommand{\N}{\mathbb{N}}
\def\Ne{\mbox{\em $\N$}}
\newcommand{\Ns}{\mbox{\tiny{\N}}}
\newcommand{\Cbar}{\bar{\mathbf{C}}}
\def\DD{\;\mbox{D}\!\!\!\!\!\!\mbox{I}\;\;\,}
\def\DDs{\mbox{{\scriptsize$\;\mbox{D}\!\!\!\!\!\!\mbox{I}\;\;\,$}}}

\def\noi{\noindent}
\def\qq{\qquad}
\def\mm{\medskip\\}
\def\lg{\langle}
\def\rg{\rangle}
\def\ra{\Righarrow}
\def\lra{\Leftrightarrow}
\def\ri{\rightarrow}
\def\fall{\mbox{ for all }}
\def\di{\displaystyle}
\def\vp{\varphi}
\def\al{\alpha}
\def\ol#1{\overline{#1}}
\def\qed{\hfill$\Box$}
\def\bref#1{(\ref{#1})}
\def\midd{\hspace*{-10pt}\left.\phantom{\di\int}\right|}
\def\text#1{\mbox{#1}}
\def\emph#1{{\em #1}}
\def\textit#1{{\em #1}}
\def\textbf#1{{\bf #1}}
\def\func#1{\mathop{\rm #1}}
\def\limfunc#1{\mathop{\rm #1}}
\def\dint{\displaystyle\int}
\def\dsum{\displaystyle\sum}
\def\dfrac{\displaystyle\frac}
\def\Bbb#1{\mathbb{#1}}
\def\bu{$\bullet$\ }
\def\sta{$\star$\ }
\def\ii{\,\mbox{i}\,}
\def\lb{\linebreak}
\def\pr{{}^{\prime}}
\def\imath{\mbox{i}}

\def\ba{\begin{array}}
\def\ea{\end{array}}
\def\beq{\begin{equation}}
\def\eeq{\end{equation}}
\def\zx#1{\begin{equation}\label{#1}}
\def\zc{\end{equation}}
\def\aru#1{\left\{{\begin{array}{l}#1\end{array}}\right.}
\def\matd#1{\left(\begin{array}{cc}#1\end{array}\right)}

\def\shw{\scriptscriptstyle}
\def\stkdn#1#2{{\mathop{#2}\limits^{}_{#1}}{}}
\def\dn#1#2{\stkdn{{\shw #1}}{#2}{}}
\def\stkup#1#2{{\mathop{#1}\limits^{#2}_{}}{}}
\def\up#1#2{\stkup{#1}{\shw #2}{}}
\def\stkud#1#2#3{{\mathop{#2}\limits^{#3}_{#1}}{}}
\def\ud#1#2#3{\stkud{\shw #1}{#2}{\shw #3}}
\def\stackunder#1#2{\mathrel{\mathop{#2}\limits_{#1}}}

\def\emptyset{\begin{picture}(13,10) \unitlength1pt
    \put(5,3){\circle{7}}\put(.7,-1.8){\line(1,1){10}}\end{picture}}

\def\fiu#1#2#3{\begin{center}
    \includegraphics[scale=#1]{#2.eps}\nopagebreak\\\parbox{10cm}
    {\begin{center}\small\bf #3\end{center}}\end{center} }
\def\fid#1#2#3#4{\begin{center}
    \includegraphics[scale=#1]{#2.eps}\hspace*{.7cm}
    \includegraphics[scale=#1]{#3.eps}
    \nopagebreak\\\parbox{14cm}{\begin{center}
    \small\bf #4\end{center}}\end{center} }
\def\fit#1#2#3#4#5{\begin{center}\begin{tabular}{ccc}
    \includegraphics[scale=#1]{#2.eps} &\hspace*{.4cm}&
    \includegraphics[scale=#1]{#3.eps}\\
    {\small\bf #4}&&{\small\bf #5}\end{tabular}\end{center}}
\pagestyle{myheadings}
\setlength{\textheight}{20cm}
\setlength{\textwidth}{13cm}
\setlength{\oddsidemargin}{18mm}
\setlength{\evensidemargin}{18mm}
\renewcommand{\theequation}{\thesection.\arabic{equation}}
\makeatletter \@addtoreset{equation}{section} \makeatother

\header{30}{36}{Dipen Ganguly, Santu Dey and Arindam Bhattacharyya}
    {On trans-Sasakian $3$-manifolds as $\eta$-Einstein solitons}
    {On trans-Sasakian $3$-manifolds as $\eta$-Einstein solitons}
\babs
    The present paper is to deliberate the class of $3$-dimensional trans-Sasakian manifolds which admits $\eta$-Einstein solitons. We have studied $\eta$-Einstein solitons on $3$-dimensional trans-Sasakian manifolds where the Ricci tensors are Codazzi type and cyclic parallel. We have also discussed some curvature conditions admitting $\eta$-Einstein solitons on $3$-dimensional trans-Sasakian manifolds and the vector field is torse-forming. We have also shown an example of $3$-dimensional trans-Sasakian manifold with respect to $\eta$-Einstein soliton to verify our results.
\eabs
\msc
    53C25, 53C15, 53C21.
\kwd
    Einstein soliton; $\eta$-Einstein soliton; Trans-Sasakian manifold; codazzi type Ricci tensor; $C$-Bochner curvature tensor.

\section{Introduction}

In 2016, G. Catino and L. Mazzieri \cite{CM} introduced the notion of Einstein soliton which can be viewed as a self-similar solution to the Einstein flow
\begin{equation}
\frac{\partial g}{\partial t}=-2(S-\frac{r}{2}g),
\end{equation}
where $g$ is the Riemannian metric, $S$ is the Ricci tensor and $r$ is the scalar curvature. It can be easily seen that the Einstein soliton is analogous to the Ricci soliton which is also generated by a self-similar solution to the very famous geometric revolution equation Ricci flow. It is a well-known fact now that the study of Ricci soliton has tremendous contribution in solving the longstanding Thurston's geometric conjecture. Similarly it is also interesting to study the Einstein soliton from various directions to solve many physical and geometrical problems. However in this paper we consider, a slight perturbation of the Einstein soliton by $\eta \otimes \eta$, called  the $\eta$-Einstein soliton. The mathematical expression for the $\eta$-Einstein soliton \cite{Blg} is given by the following equation
\begin{equation}
\mathcal{L}_{\xi}g+2S+(2\lambda-r)g+2\mu \eta \otimes \eta=0,
\end{equation}
where $\mathcal{L}_{\xi}$ denotes the Lie derivative along the direction of the vector field $\xi$, $S$ is the Ricci tensor, $r$ is the scalar curvature and $\lambda$, $\mu$ are real constants. The $\eta$-Einstein soliton is called shrinking if $\lambda <0$, steady if $\lambda =0$ and expanding if $\lambda >0$. In particular, if $\mu=0$, the $\eta$-Einstein soliton reduces to the Einstein soliton $(g,\xi,\lambda)$.
\par
 J. T. Cho and M. Kimura \cite{CK} introduced the concept of $\eta$-Ricci soliton and later C. Calin and M. Crasmareanu \cite{Cr}
studied it on Hopf hypersufaces in complex space forms. A Riemannian manifold $(M,g)$ is said to admit an $\eta$-Ricci soliton if for a smooth vector field $V$, the metric $g$ satisfies the following equation
\begin{equation}
\mathcal{L}_{\xi}g+2S+2\lambda g+2\mu \eta \otimes \eta=0,
\end{equation}
where $\mathcal{L}_{\xi}$ is the Lie derivative along the direction of $\xi$, $S$ is the Ricci tensor and $\lambda$, $\mu$ are real constants. It is to be noted that if the manifold has constant scalar curvature, then the data $(g,\xi,\lambda-\frac{r}{2},\mu)$ of the equation $(1.2)$ satisfies the equation $(1.3)$, i.e; the $\eta$-Einstein soliton reduces to an $\eta$-Ricci soliton. Hence we can remark that the two notions are different for the manifolds of non-constant scalar curvature and if the scalar curvature of the manifold is constant then the concepts of $\eta$-Ricci soliton and $\eta$-Einstein soliton coincide.
\par
\medskip

The paper is organised as follows: After a brief introduction, in section $2$, we recall some basic knowledge on trans-Sasakian manifolds. Section $3$ deals with $3$-dimensional trans-Sasakian manifolds admitting $\eta$-Einstein solitons and also the nature of the soliton is dicussed. In this section, we have constructed an example of a $3$-dimensional trans-Sasakian manifold satisfying $\eta$-Einstein soliton. In section $4$, we have contrived $\eta$-Einstein solitons in $3$-dimensional trans-Sasakian manifolds in terms of Codazzi type and cyclic parallel Ricci tensor and characterized the nature of the manifold. Sections $5,6,7,8$ are devoted to the study of some curvature conditions $R\cdot S=0$, $W_2\cdot S=0$, $R\cdot E=0$, $B\cdot S=0$, $S\cdot R=0$ admitting $\eta$-Einstein solitons in $3$-dimensional trans-Sasakian manifold. In last section we have studied torse forming vector field when $3$-dimensional trans-Sasakian manifolds admitting $\eta$-Einstein solitons.

\section{Preliminaries}
An $n$-dimensional smooth Riemannian manifold $(M,g)$ is said to be an almost contact metric manifold \cite{Bla1} if it admits a $(1,1)$ tensor field $\phi$, a characteristic vector field $\xi$, a global 1-form $\eta$ and an indefinite metric $g$ on M satisfying the following relations
\begin{eqnarray}
  \phi^2 &=& -I+\eta\otimes\xi, \\
  \eta(\xi) &=& 1, \\
  \eta(X) &=& g(X,\xi), \\
  g(\phi X,\phi Y) &=& g(X,Y)-\eta(X)\eta(Y), \\
   g(X,\phi Y)+ g(Y,\phi X) &=& 0,
\end{eqnarray}
for all vector fields $X,Y\in TM$, where $TM$ is the tangent bundle of the manifold $M$. Also it can be easily seen that $\phi(\xi)=0$, $\eta(\phi X)=0$ and rank of $\phi$ is $(n-1)$.
\par
The geometry of the almost Hermitian manifold $(M\times\mathbb{R},G,J)$ gives rise to the geometry of the almost contact metric manifold $(M,g,\phi,\xi,\eta)$, where $G$ is product metric of the product manifold $M\times\mathbb{R}$ with the complex structure $J$ defined by
\begin{equation}
  J(X,f\frac{d}{dt})=(\phi X-f\xi,\eta(X)\frac{d}{dt}),
\end{equation}
for all vector fields $X$ on the manifold $M$ and smooth function $f$ on the product manifold $M\times\mathbb{R}$. An almost contact metric manifold $(M,g,\phi,\xi,\eta)$ is called a trans-Sasakian manifold if the product manifold $(M\times\mathbb{R},G,J)$ belongs to the class $W_4$ \cite{GH}. The notion of trans-Sasakian manifolds was introduced by J. A. Oubina \cite{O} and later J. C. Marrero \cite{Ma} completely characterized the local structures of trans-Sasakian manifolds of dimension $n\geq 5$. The expression for which an almost contact metric manifold $(M,g,\phi,\xi,\eta)$ becomes a trans-Sasakian manifold is given by
\begin{equation}
  (\nabla_X\phi)(Y)=\alpha[g(X,Y)\xi-\eta(Y)X]+\beta[g(\phi X,Y)\xi-\eta(Y)\phi X],
\end{equation}
for all $X,Y\in TM$ and for some smooth functions $\alpha,\beta$ on the manifold $M$. Then such kind of manifold is called a trans-Sasakian manifold of type $(\alpha,\beta)$. In particular trans-Sasakian manifolds of type $(0,0)$, $(\alpha,0)$ and $(0,\beta)$ are called cosymplectic, $\alpha$-Sasakian and $\beta$-Kenmotsu manifolds respectively.
\par
In what follows, by a trans-Sasakian $3$-manifold, we mean a $3$-dimensional trans-Sasakian manifold $(M,g,\phi,\xi,\eta)$ of type $(\alpha,\beta)$ and we will use the notation $(M,g)$ to denote it throughout this article. Now from the above expression $(2.7)$ it can be derived that
\begin{eqnarray}
  \nabla_X\xi &=& -\alpha\phi(X)+\beta(X-\eta(X)\xi), \\
(\nabla_X\eta)(Y) &=& -\alpha g(\phi(X),Y)+\beta g(\phi(X),\phi(Y)),
\end{eqnarray}
for all vector fields $X,Y$ in $TM$. Again in a trans-Sasakian $3$-manifold $(M,g)$ the Ricci tensor is given by
\begin{eqnarray}
  S(X,Y) &=& [\frac{r}{2}+\xi\beta-(\alpha^2-\beta^2)]g(X,Y)-[\frac{r}{2}+\xi\beta-3(\alpha^2-\beta^2)]\eta(X)\eta(Y), \nonumber\\
   && -[Y\beta+(\phi(Y)\alpha)]\eta(X)-[X\beta+(\phi(X)\alpha)]\eta(Y).
\end{eqnarray}
Furthermore, if the functions $\alpha$, $\beta$ are constants then, in a trans-Sasakian $3$-manifold $(M,g)$ the following relations hold,
\begin{eqnarray}
  R(X,Y)\xi &=& (\alpha^2-\beta^2)[\eta(Y)X-\eta(X)Y], \\
  R(\xi,X)Y &=& (\alpha^2-\beta^2)[g(X,Y)\xi-\eta(Y)X], \\
  R(\xi,X)\xi &=& (\alpha^2-\beta^2)[\eta(X)\xi-X], \\
  S(X,Y) &=& [\frac{r}{2}-(\alpha^2-\beta^2)]g(X,Y)-[\frac{r}{2}-3(\alpha^2-\beta^2)]\eta(X)\eta(Y), \\
  S(X,\xi) &=& 2(\alpha^2-\beta^2)\eta(X),
\end{eqnarray}
for all vector fields $X,Y$ in $TM$ and where $R$ is the curvature tensor and $S$ is the Ricci tensor.

\begin{definition}
  A trans-Sasakian $3$-manifold $(M,g)$ is said to be an $\eta$-Einstein manifold if its Ricci tensor $S$ is of the form
  \begin{equation}
  S(X,Y)=ag(X,Y)+b\eta(X)\eta(Y),
  \end{equation}
  for all $X,Y\in TM$ and smooth functions $a,b$ on the manifold $(M,g)$.
\end{definition}

\section{$\eta$-Einstein solitons on trans-Sasakian $3$-manifolds}
Let us consider a trans-Sasakian $3$-manifold $(M,g)$ admitting an $\eta$-Einstein soliton $(g,\xi,\lambda,\mu)$. Then from equation $(1.2)$ we can write
\begin{equation}
(\mathcal{L}_{\xi}g)(X,Y)+2S(X,Y)+[2\lambda -r]g(X,Y)+2\mu \eta(X)\eta(Y)=0,
\end{equation}
for all $X,Y\in TM$.\\
Again from the well-known formula $(\mathcal{L}_{\xi}g)(X,Y)=g(\nabla_X\xi,Y)+g(\nabla_Y\xi,X)$ of Lie-derivative and using $(2.8)$, we obtain for a trans-Sasakian $3$-manifold
\begin{equation}
(\mathcal{L}_{\xi}g)(X,Y)=2\beta[g(X,Y)-2\beta\eta(X)\eta(Y)].
\end{equation}
Now in view of the equations $(3.1)$ and $(3.2)$ we get
\begin{equation}
S(X,Y)=(\frac{r}{2}-\lambda-\beta)g(X,Y)+(\beta-\mu)\eta(X)\eta(Y).
\end{equation}
This shows that the manifold $(M,g)$ is an $\eta$-Einstein manifold.\\
Also from equation $(3.3)$ replacing $Y=\xi$ we find that
\begin{equation}
S(X,\xi)=(\frac{r}{2}-\lambda-\mu)\eta(X).
\end{equation}
Compairing the above equation $(3.4)$ with $(2.15)$ yields
\begin{equation}
r=4(\alpha^2-\beta^2)+2\lambda+2\mu.
\end{equation}
Again, considering an orthonormal basis $\{e_1,e_2,e_3\}$ of $(M,g)$ and then setting $X=Y=e_i$ in equation $(3.3)$ and summing over $i=1,2,3$ we get
\begin{equation}
r=6\lambda+4\beta+2\mu.
\end{equation}
Finally combining equations $(3.5)$ and $(3.6)$ we arrive at
\begin{equation}
\lambda=(\alpha^2-\beta^2)-\beta.
\end{equation}
Thus the above discussion leads to the following
\begin{theorem}
  If a trans-Sasakian $3$-manifold $(M,g)$ admits an $\eta$-Einstein soliton $(g,\xi,\lambda,\mu)$, then the manifold $(M,g)$ becomes an $\eta$-Einstein manifold of constant scalar curvature $r=6\lambda+4\beta+2\mu$. Furthermore, the soliton is shrinking, steady or expanding according as; $\alpha^2<\beta(\beta+1)$, $\alpha^2=\beta(\beta+1)$, $\alpha^2>\beta(\beta+1)$ respectively.
\end{theorem}
\medskip

\begin{example}
Let us consider the $3$-dimensional manifold $M=\{(u,v,w)\in\mathbb{R}^3:w\neq 0\}$. Define a linearly independent set of vector fields
$\{e_i: 1\leq i\leq 3\}$ on the manifold $M$ given by
\begin{equation}
e_1=e^{2w}\frac{\partial}{\partial u}, ~~~~e_2=e^{2w}\frac{\partial}{\partial v}, ~~~~e_3=\frac{\partial}{\partial w}.\nonumber
\end{equation}
Let us define the Riemannian metric $g$ on $M$ by
\begin{equation}
g(e_i,e_j)=\left\{ \begin{array}{rcl}
1, & \mbox{for}
& i=j \\ 0, & \mbox{for} & i\neq j
\end{array}\right.\nonumber
\end{equation}
for all $i,j=1,2,3$. Now considering $e_3=\xi$, let us take the $1$-form $\eta$, on the manifold $M$, defined by
\begin{equation}
\eta(U)=g(U,e_3),~~~\forall U\in TM.\nonumber
\end{equation}
Then it can be observed that $\eta(\xi)=1$. Let us define the $(1,1)$ tensor field $\phi$ on $M$ as
\begin{equation}
\phi(e_1)=e_2, ~~~~\phi(e_2)=-e_1, ~~~~\phi(e_3)=0.\nonumber
\end{equation}
Using the linearity of $g$ and $\phi$ it can be easily checked that
\begin{equation}
\phi^2(U)=-U+\eta(U)\xi,~~g(\phi U,\phi V)=g(U,V)-\eta(U)\eta(V),~~~\forall U,V\in TM.\nonumber
\end{equation}
Hence the structure $(g,\phi, \xi, \eta)$ defines an almost contact metric structure on the manifold $M$. Now, using the definitions of Lie bracket, after some direct computations we get
\begin{equation}
[e_1,e_2]=0,~~~~[e_1,e_3]=-2e_1,~~~~[e_2,e_3]=-2e_2.\nonumber
\end{equation}
Again the Riemannian connection $\nabla$ of the metric $g$ is defined by the well-known Koszul's formula which is given by
\begin{eqnarray}
  2g(\nabla_XY,Z) &=& Xg(Y,Z)+Yg(Z,X)-Zg(X,Y) \nonumber\\
   && -g(X,[Y,Z])+g(Y,[Z,X])+g(Z,[X,Y]).\nonumber
\end{eqnarray}
Using the above formula one can easily calculate that
\begin{multline}
 ~~~~~~~~~~~~~~~~~~~~~~~~\nabla_{e_1}e_1=2e_3,~~~\nabla_{e_1}e_2=0,~~~\nabla_{e_1}e_3=-2e_1, \\
 \nabla_{e_2}e_1=0,~~~\nabla_{e_2}e_2=2e_3,~~~\nabla_{e_2}e_3=-2e_2, \\
 \nabla_{e_3}e_1=0,~~~\nabla_{e_3}e_2=0,~~~\nabla_{e_3}e_3=0.~~~~~~~~~~~~~~~~~~~~~~~~~~~~~~~~\nonumber
\end{multline}
 Thus from the above relations it follows that the manifold $(M,g)$ is a trans-Sasakian $3$-manifold. Now using the well-known formula $R(X,Y)Z=\nabla_X\nabla_YZ-\nabla_Y\nabla_XZ-\nabla_{[X,Y]}Z$ the non-vanishing components of the Riemannian curvature tensor $R$ can be easily obtained as
\begin{multline}
  ~~~~~~~~~~~~~~~~~~~~~~~~~~~~~~R(e_1,e_2)e_2=R(e_1,e_3)e_3=-4e_1, \\
  R(e_2,e_3)e_3=R(e_3,e_1)e_1=-4e_2, \\
  R(e_3,e_2)e_2=4e_2,~~~~R(e_2,e_1)e_1=4e_3.~~~~~~~~~~~~~~~~~~~~~~~~~~\nonumber
\end{multline}
Hence we can calculate the components of the Ricci tensor as follows
\begin{equation}
S(e_1,e_1)=0,~~~S(e_2,e_2)=0,~~~S(e_3,e_3)=-8.\nonumber
\end{equation}
Therefore in view of the above values of the Ricci tensor, from the equation $(1.2)$  we can calculate $\lambda=-2$ and $\mu=6$. Hence we can say that the data $(g,\xi,-2,6)$ defines an $\eta$-Einstein soliton on the trans-Sasakian $3$-manifold $(M,g)$. Also we can see that the manifold $(M,g)$ is a manifold of constant scalar curvature $r=-8$ and hence the theorem $(3.1)$ is verified.
\end{example}
\par
\medskip

Next we consider a trans-Sasakian $3$-manifold $(M,g)$ and assume that it admits an $\eta$-Einstein soliton $(g,V,\lambda,\mu)$ such that $V$ is pointwise collinear with $\xi$, i.e; $V=b\xi$, for some function $b$; then from the equation $(1.2)$ it follows that
\begin{multline}
  ~~~~~~~~~~~bg(\nabla_X\xi,Y)+(Xb)\eta(Y)+bg(\nabla_Y\xi,X)+(Yb)\eta(X) \nonumber\\
  +2S(X,Y)+(2\lambda -r)g(X,Y)+2\mu \eta(X)\eta(Y)=0.~~~~~~~~~~~~~~
\end{multline}
Then using the equation $(2.8)$ in above we arrive at
\begin{multline}
  ~~~~~~~~~~~(2b\beta +2\lambda -r)g(X,Y)+(Xb)\eta(Y)+(Yb)\eta(X)\\
  +2S(X,Y)+2(b\beta+\mu)\eta(X)\eta(Y)=0.~~~~~~~~~~~~~~
\end{multline}
Replacing $Y=\xi$ in the above equation yields
\begin{equation}
(Xb)+(\xi b)\eta(X)+2S(X,\xi)+(2\lambda +2\mu-r)\eta(X)=0.
\end{equation}
Again taking $X=\xi$ in $(3.9)$ and by virtue of $(2.15)$ we arrive at
\begin{equation}
2(\xi b)=(r-2\lambda-2\mu)-4(\alpha^2-\beta^2).
\end{equation}
Using this value from $(3.10)$ in the equation $(3.9)$ and recalling $(2.15)$ we can write
\begin{equation}
db=[\frac{r}{2}-\lambda-\mu-2(\alpha^2-\beta^2)]\eta.
\end{equation}
Now taking exterior differentiation on both sides of $(3.11)$ and using the famous Poincare's lemma i.e; $d^2=0$, finally we arrive at
\begin{equation}
r=2\lambda+2\mu+4(\alpha^2-\beta^2).
\end{equation}
In view of the above $(3.12)$ the equation $(3.11)$ gives us $db=0$ i.e; the function $b$ is constant. Then the equation $(3.8)$ reduces to
\begin{equation}
S(X,Y)=(\frac{r}{2}-\lambda-b\beta)g(X,Y)+(b\beta-\mu)\eta(X)\eta(Y),
\end{equation}
for all $X,Y\in TM$. Hence we can state the following
\begin{theorem}
   If a trans-Sasakian $3$-manifold $(M,g)$ admits an $\eta$-Einstein soliton $(g,V,\lambda,\mu)$ such that $V$ is pointwise collinear with $\xi$, then $V$ is constant multiple of $\xi$ and the manifold $(M,g)$ becomes an $\eta$-Einstein manifold of constant scalar curvature $r=2\lambda+2\mu+4(\alpha^2-\beta^2)$.
\end{theorem}

\section{$\eta$-Einstein solitons on trans-Sasakian $3$-manifolds with Codazzi type and cyclic parallel Ricci tensor}
The purpose of this section is to study $\eta$-Einstein solitons in trans-Sasakian $3$-manifolds having certain special types of Ricci tensor namely codazzi type Ricci tensor and cyclic parallel Ricci tensor.
\begin{definition}
  \cite{G} A trans-Sasakian $3$-manifold is said to have Codazzi type Ricci tensor if its Ricci tensor $S$ is non-zero and satisfies the following relation
  \begin{equation}
  (\nabla_XS)(Y,Z)=(\nabla_YS)(X,Z),~~~\forall X,Y,Z\in TM.
  \end{equation}
\end{definition}
Let us consider a trans-Sasakian $3$-manifold having Codazzi type Ricci tensor and admits an $\eta$-Einstein soliton $(g,\xi,\lambda,\mu)$, then equation $(3.3)$ holds. Now covariantly differentiating the equation $(3.3)$ with respect to an arbitrary vector field $X$ and then using $(2.9)$ we get
\begin{eqnarray}
  (\nabla_XS)(Y,Z) &=& 2(\beta-\mu)[\eta(Y)(-\alpha g(\phi X,Z)+\beta g(\phi X,\phi Z)) \nonumber\\
   && +\eta(Z)(-\alpha g(\phi X,Y)+\beta g(\phi X,\phi Y))].
\end{eqnarray}
  Similarly we can compute
  \begin{eqnarray}
  (\nabla_YS)(X,Z) &=& 2(\beta-\mu)[\eta(X)(-\alpha g(\phi Y,Z)+\beta g(\phi Y,\phi Z)) \nonumber\\
   && +\eta(Z)(-\alpha g(\phi Y,X)+\beta g(\phi Y,\phi X))].
\end{eqnarray}
Since the manifold has Codazzi type Ricci tensor, using $(4.2)$ and $(4.3)$ in the equation $(4.1)$ and then recalling $(2.4)$ we arrive at
\begin{eqnarray}
  2(\beta-\mu)[\eta(Y)(-\alpha g(\phi X,Z)+\beta g(X,Z))-\eta(X)(-\alpha g(\phi Y,Z) &&  \nonumber\\
  +\beta g(Y,Z))-2\alpha\eta(Z)g(\phi X,Y)] &=& 0.
\end{eqnarray}
Putting $Z=\xi$ in above and in view of $(2.2)$ we finally obtain
\begin{equation}
4\alpha(\mu-\beta)g(\phi X,Y)=0,
\end{equation}
forall $X,Y\in TM$. Therefore from $(4.4)$ we can conclude that either $\alpha=0$ or $\mu=\beta$. Hence we have the following
\begin{theorem}
  Let $(M,g)$ be  a trans-Sasakian $3$-manifold admitting an $\eta$-Einstein soliton $(g,\xi,\lambda,\mu)$. If the Ricci tensor of the manifold is of Codazzi type then the manifold becomes a $\beta$-Kenmotsu manifold provided $\mu\neq\beta$.
\end{theorem}
Now using $\alpha=0$ in equation $(3.7)$ we get $\lambda=-\beta(\beta+1)$. Thus we can state the following
\begin{corollary}
 Let $(M,g)$ be  a trans-Sasakian $3$-manifold admitting an $\eta$-Einstein soliton $(g,\xi,\lambda,\mu)$ with $\mu\neq\beta$. If the Ricci tensor of the manifold is of Codazzi type then the soliton is shrinking if $\beta<-1$ or, $\beta>0$; steady if $\beta=-1$ or $\beta=0$; and expanding if $-1<\beta<0$ respectively.
\end{corollary}
Again from the equation $(4.4)$ we can write that $\mu=\beta$ if $\alpha\neq 0$. Then from equation $(3.3)$ we obtain
\begin{equation}
S(X,Y)=(\frac{r}{2}-\lambda-\beta)g(X,Y),
\end{equation}
for all $X,Y\in TM$. Then contracting the equation $(4.5)$ we get $r=6\lambda+6\beta$. Hence in view of this and equation $(4.5)$ we have the following
\begin{theorem}
  Let $(M,g)$ be a trans-Sasakian $3$-manifold admitting an $\eta$-Einstein soliton $(g,\xi,\lambda,\mu)$. If the Ricci tensor of the manifold is of Codazzi type then the manifold becomes an Einstein manifold of constant scalar curvature $r=6\lambda+6\beta$ provided $\alpha\neq 0$.
\end{theorem}
\medskip
\begin{definition}
  \cite{G} A trans-Sasakian $3$-manifold is said to have cyclic parallel Ricci tensor if its Ricci tensor $S$ is non-zero and satisfies the following relation
  \begin{equation}
  (\nabla_XS)(Y,Z)+(\nabla_YS)(Z,X)+(\nabla_ZS)(X,Y)=0~~~\forall X,Y,Z\in TM.
  \end{equation}
\end{definition}
\par
\medskip
Let us consider a trans-Sasakian $3$-manifold, having cyclic parallel Ricci tensor, admits an $\eta$-Einstein soliton $(g,\xi,\lambda,\mu)$, then equation $(3.3)$ holds. Now taking covariant differentiation of $(3.3)$ and using equation $(2.9)$ we obtain relations $(4.2)$ and $(4.3)$. In a similar manner we get the following
\begin{eqnarray}
  (\nabla_ZS)(X,Y) &=& 2(\beta-\mu)[\eta(X)(-\alpha g(\phi Z,Y)+\beta g(\phi Z,\phi Y)) \nonumber\\
   && +\eta(Y)(-\alpha g(\phi Z,X)+\beta g(\phi Z,\phi X))].
\end{eqnarray}
Now since the manifold has cyclic parallel Ricci tensor, using the values from $(4.2)$, $(4.3)$ and $(4.8)$ in the equation $(4.7)$ and then making use of $(2.4)$ we arrive at
\begin{equation}
  4\beta(\beta-\mu)[\eta(X)g(\phi Y,\phi Z)+\eta(Y)g(\phi Z,\phi X)+\eta(Z)g(\phi X,\phi Y)]=0.
  \end{equation}
  Replacing $Z=\xi$ in the above equation $(4.9)$ yields
  \begin{equation}
 4\beta(\beta-\mu)g(\phi X,\phi Y)=0,
  \end{equation}
  for all $X,Y\in TM$. Since $g(\phi X,\phi Y)\neq 0$ always, the above equation $(4.10)$ implies that either $\beta =0$ or, $\mu=\beta$. Thus we can state the following
  \begin{theorem}
  Let $(M,g)$ be  a trans-Sasakian $3$-manifold admitting an $\eta$-Einstein soliton $(g,\xi,\lambda,\mu)$. If the manifold has cyclic parallel Ricci tensor, then the manifold becomes an $\alpha$-Sasakian manifold provided $\mu\neq\beta$.
\end{theorem}
Now using $\beta=0$ in equation $(3.7)$ we get $\lambda=\alpha^2>0$. Therefore we have
\begin{corollary}
 Let $(M,g)$ be  a trans-Sasakian $3$-manifold admitting an $\eta$-Einstein soliton $(g,\xi,\lambda,\mu)$ with $\mu\neq\beta$. If the manifold has cyclic parallel Ricci tensor then the soliton is expanding.
\end{corollary}
Again if $\beta\neq 0$ then from $(4.10)$ it follows that $\mu=\beta$. Therefore after a similar calculation like equation $(4.6)$ we can state
\begin{theorem}
  Let $(M,g)$ be  a trans-Sasakian $3$-manifold admitting an $\eta$-Einstein soliton $(g,\xi,\lambda,\mu)$. If the manifold has cyclic parallel Ricci tensor, then the manifold becomes an Einstein manifold of constant scalar curvature $r=6\lambda+6\beta$ provided $\beta\neq 0$.
\end{theorem}

\section{$\eta$-Einstein solitons on trans-Sasakian $3$-manifolds satisfying $R(\xi,X)\cdot S=0$ and $W_2(\xi,X)\cdot S=0$}

Let us first consider a trans-Sasakian $3$-manifold which admits an $\eta$-Einstein soliton $(g,\xi,\lambda,\mu)$ and the manifold satisfies the curvature condition $R(\xi,X)\cdot S=0$. Then $\forall X,Y,Z\in TM$ we can write
\begin{equation}
S(R(\xi,X)Y,Z)+S(Y,R(\xi,X)Z)=0.
\end{equation}
Now using $(3.3)$ in $(5.1)$ we get
\begin{eqnarray}
  (\frac{r}{2}-\lambda-\beta)g(R(\xi,X)Y,Z)+(\beta-\mu)\eta(R(\xi,X)Y)\eta(Z)~~~~~~~~~ &&  \nonumber\\
 ~~~~~~~~~~~+(\frac{r}{2}-\lambda-\beta)g(R(\xi,X)Z,Y)+(\beta-\mu)\eta(R(\xi,X)Z)\eta(Y)=0. &&
\end{eqnarray}
In view of $(2.12)$ the previous equation becomes
\begin{equation}
(\alpha^2-\beta^2)(\beta-\mu)[g(X,Y)\eta(Z)+g(X,Z)\eta(Y)-2\eta(X)\eta(Y)\eta(Z)]=0.
\end{equation}
Putting $Z=\xi$ in the above equation $(5.3)$ and recalling $(2.4)$ obtain
\begin{equation}
(\alpha^2-\beta^2)(\beta-\mu)g(\phi X,\phi Y)=0,
\end{equation}
for all $X,Y\in TM$. Since $g(\phi X,\phi X)\neq 0$ always and for non-trivial case $\alpha^2\neq\beta^2$, we can conclude from the equation $(5.4)$ that $\mu=\beta$. Then from equation $(3.3)$ we obtain
\begin{equation}
S(X,Y)=(\frac{r}{2}-\lambda-\beta)g(X,Y),
\end{equation}
for all $X,Y\in TM$. Then contracting the equation $(5.5)$ we get $r=6\lambda+6\beta$. Hence in view of this and equation $(5.5)$ we have the following
\begin{theorem}
  Let $(M,g)$ be a trans-Sasakian $3$-manifold admitting an $\eta$-Einstein soliton $(g,\xi,\lambda,\mu)$. If the manifold satisfies the curvature condition $R(\xi,X)\cdot S=0$, then the manifold becomes an Einstein manifold of constant scalar curvature $r=6\lambda+6\beta$.
\end{theorem}
\par
\medskip
Our next result of this section is on $W_2$-curvature tensor. It is an important curvature tensor which was introduced in 1970 by Pokhariyal and Mishra \cite{PM}. For this let us recall the definition of $W_2$-curvature tensor as follows
\begin{definition}
  The $W_2$-curvature tensor in a trans-Sasakian $3$-manifold $(M,g)$ is defined as
\begin{equation}
  W_2(X,Y)Z=R(X,Y)Z+\frac{1}{2}[g(X,Z)QY-g(Y,Z)QX].
  \end{equation}
\end{definition}
Now assume that $(M,g)$ is a trans-Sasakian $3$-manifold admitting an $\eta$-Einstein soliton $(g,\xi,\lambda,\mu)$ and also the manifold satisfies the curvature condition $W_2(\xi,X)\cdot S=0$. Then we can write
\begin{equation}
S(W_2(\xi,X)Y,Z)+S(Y,W_2(\xi,X)Z)=0,~~~\forall X,Y,Z\in TM.
\end{equation}
In view of $(3.3)$ the above equation $(5.7)$ becomes
\begin{eqnarray}
  (\frac{r}{2}-\lambda-\beta)[g(W_2(\xi,X)Y,Z)+g(W_2(\xi,X)Z,Y)]~~~~~~~~~ &&  \nonumber\\
 ~~~~~~~~~~~+(\beta-\mu)[\eta(W_2(\xi,X)Y)\eta(Z)+\eta(W_2(\xi,X)Z)\eta(Y)]=0. &&
\end{eqnarray}
Again from $(3.3)$ it follows that
\begin{equation}
QX=(\frac{r}{2}-\lambda-\beta)X+(\beta-\mu)\eta(X)\xi,
\end{equation}
which implies
\begin{equation}
Q\xi=(\frac{r}{2}-\lambda-\mu)\xi.
\end{equation}
Replacing $X=\xi$ in $(5.6)$ and then using equations $(2.12)$, $(5.9)$ and $(5.10)$ we obtain
\begin{equation}
W_2(\xi,Y)Z=Bg(Y,Z)\xi-A\eta(Z)Y+(A-B)\eta(Y)\eta(Z),
\end{equation}
where $A=(\alpha^2-\beta^2)-\frac{1}{2}(\frac{r}{2}-\lambda-\beta)$ and $B=(\alpha^2-\beta^2)-\frac{1}{2}(\frac{r}{2}-\lambda-\mu)$. Taking inner product of $(5.11)$ with respect to the vector field $\xi$ yields
\begin{equation}
\eta(W_2(\xi,Y)Z)=B[g(Y,Z)-\eta(Y)\eta(Z)].
\end{equation}
Using $(5.11)$ and $(5.12)$ in the equation $(5.8)$ and then taking $Z=\xi$ we arrive at
\begin{equation}
(A-B)[2B-(\frac{r}{2}-\lambda-\beta)][g(X,Y)-\eta(X)\eta(Y)]=0,\nonumber
\end{equation}
which in view of $(2.4)$ implies
\begin{equation}
(A-B)[2B-(\frac{r}{2}-\lambda-\beta)]g(\phi X,\phi Y)=0,
\end{equation}
for all $X,Y\in TM$. Since $g(\phi X,\phi X)\neq 0$ always, we can conclude from the equation $(5.13)$ that either $A=B$ or, $2B=\frac{r}{2}-\lambda-\beta$. Thus recalling the values of $A$ and $B$ it implies that either $\mu=\beta$ or,
\begin{equation}
2(\alpha^2-\beta^2)=r-2\lambda-\mu-\beta.
\end{equation}
Now for the case $\mu=\beta$, proceeding similarly as the equation $(5.5)$ we can say that the manifold becomes an Einstein manifold. Again combining $(5.14)$ with $(3.5)$ we get
\begin{equation}
r=2\lambda+2\beta.
\end{equation}
Therefore we can state the following
\begin{theorem}
  Let $(M,g)$ be a trans-Sasakian $3$-manifold admitting an $\eta$-Einstein soliton $(g,\xi,\lambda,\mu)$. If the manifold satisfies the curvature condition $W_2(\xi,X)\cdot S=0$, then either the manifold becomes an Einstein manifold or it is a manifold of constant scalar curvature $r=2\lambda+2\beta$.
\end{theorem}
Again in view of $(3.6)$, the equation $(5.15)$ implies $\lambda=-\frac{1}{2}(\mu+\beta)$. Hence we have
\begin{corollary}
  Let $(M,g)$ be a trans-Sasakian $3$-manifold admitting an $\eta$-Einstein soliton $(g,\xi,\lambda,\mu)$ with $\mu\neq\beta$. If the manifold satisfies the curvature condition $W_2(\xi,X)\cdot S=0$, then the soliton is expanding, steady or shrinking according as $\mu<-\beta$, $\mu=-\beta$ or, $\mu>-\beta$ respectively.
\end{corollary}

\section{Einstein semi-symmetric trans-Sasakian $3$-manifolds admitting $\eta$-Einstein solitons}
\begin{definition}
  A trans-Sasakian $3$-manifold $(M,g)$ is called Einstein semi-symmetric \cite{Sz} if $R.E=0$, where $E$ is the Einstein tensor given by
  \begin{equation}
  E(X,Y)=S(X,Y)-\frac{r}{3}g(X,Y),
  \end{equation}
  for all vector fields $X,Y\in TM$ and $r$ is the scalar curvature of the manifold.
\end{definition}
Now consider a trans-Sasakian $3$-manifold is Einstein semi-symmetric i.e; the manifold satisfies the curvature condition $R.E=0$. Then for all vector fields $X,Y,Z,W\in TM$ we can write
\begin{equation}
E(R(X,Y)Z,W)+E(Z,R(X,Y)W)=0.
\end{equation}
In view of $(6.1)$ the equation $(6.2)$ becomes
\begin{equation}
S(R(X,Y)Z,W)+S(Z,R(X,Y)W)=\frac{r}{3}[g(R(X,Y)Z,W)+g(Z,R(X,Y)W)].
\end{equation}
Replacing $X=Z=\xi$ in the above equation $(6.3)$ and then using $(2.12)$, $(2.13)$ we arrive at
\begin{equation}
(\alpha^2-\beta^2)S(Y,W)=(\alpha^2-\beta^2)[\eta(Y)S(\xi,W)+\eta(W)S(\xi,Y)-g(Y,W)S(\xi,\xi)].
\end{equation}
So, now in view of $(2.15)$ the above equation $(6.4)$ finally yields
\begin{equation}
S(Y,W)=-2(\alpha^2-\beta^2)g(Y,W)+4(\alpha^2-\beta^2)\eta(Y)\eta(W),
\end{equation}
for all $Y,W\in TM$. This implies that the manifold is an $\eta$-Einstein manifold. Hence we have the following
\begin{lemma}
  An Einstein semi-symmetric trans-Sasakian $3$-manifold is an $\eta$-Einstein manifold.
\end{lemma}
Now let us assume that the Einstein semi-symmetric trans-Sasakian $3$-manifold $(M,g)$ admits an $\eta$-Einstein soliton $(g,\xi,\lambda,\mu)$. Then equation $(3.3)$ holds and combining $(3.3)$ with the above equation $(6.5)$ we get
\begin{equation}
r=2\lambda+\mu+\beta.
\end{equation}
Again recalling the equation $(3.6)$ in the above $(6.6)$ we have
\begin{equation}
\lambda=-\frac{1}{4}(\mu+3\beta).
\end{equation}
Therefore we can state the following
\begin{theorem}
  Let $(M,g)$ be a trans-Sasakian $3$-manifold admitting an $\eta$-Einstein soliton $(g,\xi,\lambda,\mu)$. If the manifold is Einstein semi-symmetric, then the manifold becomes an $\eta$-Einstein manifold of constant scalar curvature $r=2\lambda+\mu+\beta$ and the soliton is expanding, steady or shrinking according as $\mu<3\beta$, $\mu=3\beta$ or, $\mu>3\beta$ respectively.
\end{theorem}

\section{$\eta$-Einstein solitons on trans-Sasakian $3$-manifolds satisfying $B(\xi,X)\cdot S=0$}

In 1949, S. Bochner \cite{Boc} introduced the concept of the well-known Bochner curvature tensor merely as a K\"{a}hler analogue of the Weyl conformal curvature tensor but the geometric significance of it in the light of Boothby-Wangs fibration was presented later by D. E. Blair \cite{Bla}. The notion of C-Bochner curvature tensor in a Sasakian manifold was introduced by M. Matsumoto, G. Ch\={u}man \cite{MC} in 1969. The C-Bochner curvature tensor in trans-Sasakian $3$-manifold $(M,g)$ is given by
\begin{eqnarray}
  B(X,Y)Z &=& R(X,Y)Z+\frac{1}{6}[g(X,Z)QY-S(Y,Z)-g(Y,Z)QX \nonumber\\
   && +S(X,Z)Y+g(\phi X,Z)Q\phi Y-S(\phi Y,Z)\phi X-g(\phi Y,Z)Q\phi X \nonumber\\
   && +S(\phi X,Z)\phi Y+2S(\phi X,Y)\phi Z+2g(\phi X,Y))Q\phi Z \nonumber\\
   && +\eta(Y)\eta(Z)QX-\eta(Y)S(X,Z)\xi+\eta(X)S(Y,Z)\xi-\eta(X)\eta(Z)QY] \nonumber\\
   && -\frac{D+2}{6}[g(\phi X,Z)\phi Y-g(\phi Y,Z)\phi X+2g(\phi X,Y)\phi Z] \nonumber\\
   && +\frac{D}{6}[\eta(Y)g(X,Z)\xi-\eta(Y)\eta(Z)X+\eta(X)\eta(Z)Y-\eta(X)g(Y,Z)\xi] \nonumber\\
   && -\frac{D-4}{6}[g(X,Z)Y-g(Y,Z)X],
\end{eqnarray}
where $D=\frac{r+2}{4}$.\par
\medskip
Let us consider a trans-Sasakian $3$-manifold $(M,g)$ which admits an $\eta$-Einstein soliton $(g,\xi,\lambda,\mu)$ and also the manifold satisfies the curvature condition $B(\xi,X)\cdot S=0$. Then $\forall X,Y,Z\in TM$ we can write
\begin{equation}
S(B(\xi,X)Y,Z)+S(Y,B(\xi,X)Z)=0.
\end{equation}
Now using $(3.3)$ in $(7.2)$ we get
\begin{eqnarray}
  (\frac{r}{2}-\lambda-\beta)[g(B(\xi,X)Y,Z)+g(B(\xi,X)Z,Y)]~~~~~~~~~ &&  \nonumber\\
 ~~~~~~~~~~~+(\beta-\mu)[\eta(B(\xi,X)Y)\eta(Z)+\eta(B(\xi,X)Z)\eta(Y)]=0. &&
\end{eqnarray}
Again from $(3.3)$ it follows that
\begin{equation}
QX=(\frac{r}{2}-\lambda-\beta)X+(\beta-\mu)\eta(X)\xi,
\end{equation}
which implies
\begin{equation}
Q\xi=(\frac{r}{2}-\lambda-\mu)\xi.
\end{equation}
Also taking $X=\xi$ in $(7.1)$ we obtain
\begin{eqnarray}
  B(\xi,Y)Z &=& R(\xi,Y)Z\frac{1}{6}[S(\xi,Z)Y-g(Y,Z)Q\xi+\eta(Y)\eta(Z)Q\xi \nonumber\\
   && -\eta(Y)S(\xi,Z)\xi]+\frac{4}{6}[\eta(Z)Y-g(Y,Z)\xi].
\end{eqnarray}
Using equations $(2.12)$, $(3.4)$ and $(7.5)$ in the above equation $(7.6)$ yields
\begin{equation}
B(\xi,Y)Z=[(\alpha^2-\beta^2)-\frac{1}{6}(\frac{r}{2}-\lambda-\mu)-\frac{4}{6}][g(Y,Z)\xi-\eta(Z)Y].
\end{equation}
In view of $(7.7)$ the equation $(7.3)$ becomes
\begin{eqnarray}
[(\alpha^2-\beta^2)-\frac{1}{6}(\frac{r}{2}-\lambda-\mu)-\frac{4}{6}](\beta-\mu)[g(X,Y)\eta(Z) && \nonumber\\
 +g(X,Z)\eta(Y)-2\eta(X)\eta(Y)\eta(Z)] &=& 0. \nonumber
\end{eqnarray}
Replacing $Z=\xi$ in the above equation and recalling $(2.4)$, finally we arrive at
\begin{equation}
[(\alpha^2-\beta^2)-\frac{1}{6}(\frac{r}{2}-\lambda-\mu)-\frac{4}{6}](\beta-\mu)g(\phi X,\phi Y)=0,
\end{equation}
for all vector fields $X,Y\in TM$. Hence from $(7.8)$ we can conclude that either
\begin{equation}
[(\alpha^2-\beta^2)-\frac{1}{6}(\frac{r}{2}-\lambda-\mu)-\frac{4}{6}]=0,
\end{equation}
or, $\mu=\beta$. Also for $\mu=\beta$ proceeding similarly as equation $(4.6)$ it can be easily shown that the manifold becomes an Einstein manifold. Again if $\mu\neq\beta$ using $(3.7)$ in the equation $(7.9)$ we have
\begin{equation}
r=10\lambda+2\mu+12\beta-8,
\end{equation}
which implies that the manifold becomes a manifold of constant scalar curvature. Therefore we can state the following
\begin{theorem}
  Let $(M,g)$ be a trans-Sasakian $3$-manifold admitting an $\eta$-Einstein soliton $(g,\xi,\lambda,\mu)$. If the manifold satisfies the curvature condition $B(\xi,X)\cdot S=0$, then either the manifold is an Einstein manifold or it is a  manifold of constant scalar curvature $r=10\lambda+2\mu+12\beta-8$.
\end{theorem}
Now for the case $\mu\neq\beta$, using the equation $(3.6)$ in $(7.10)$ we obtain $\lambda=2(1-\beta)$. Hence we have
\begin{corollary}
  Let $(M,g)$ be a trans-Sasakian $3$-manifold admitting an $\eta$-Einstein soliton $(g,\xi,\lambda,\mu)$ with $\mu\neq\beta$. If the manifold satisfies the curvature condition $B(\xi,X)\cdot S=0$, then the soliton is expanding, steady or shrinking according as $\beta<1$, $\beta=1$ or, $\beta>1$ respectively.
\end{corollary}

\section{$\eta$-Einstein solitons on trans-Sasakian $3$-manifolds satisfying $S(\xi,X)\cdot R=0$}
In this section we study the curvature condition $S(\xi,X)\cdot R=0$, where by $\cdot$ we denote the derivation of the tensor algebra at each point of the tangent space as follows:
\begin{eqnarray}
  S((\xi,X)\cdot R)(Y,Z)W &:=& ((\xi\wedge_SX)\cdot R)(Y,Z)W \nonumber\\
   &:=& (\xi\wedge_SX)R(Y,Z)W+R((\xi\wedge_SX)Y,Z)W \nonumber\\
   && +R(Y,(\xi\wedge_SX)Z)W+R(Y,Z)(\xi\wedge_SX)W,
\end{eqnarray}
where the endomorphism $X\wedge_SY$ is defined by
\begin{equation}
(X\wedge_SY)Z:=S(Y,Z)X-S(X,Z)Y.\nonumber
\end{equation}
Now let us consider a trans-Sasakian $3$-manifold $(M,g)$ which admits an $\eta$-Einstein soliton $(g,\xi,\lambda,\mu)$ and also the manifold satisfies the curvature condition $S(\xi,X)\cdot R=0$. Then using this condition and the equation $(8.1)$ we can write
\begin{eqnarray}
  S(X,R(Y,Z)W)\xi-S(\xi,R(Y,Z)W)X+S(X,Y)R(\xi,Z)W &&  \nonumber\\
  -S(\xi,Y)R(X,Z)W+S(X,Z)R(Y,\xi)W-S(\xi,Z)R(Y,X)W &&  \nonumber\\
  +S(X,W)R(Y,Z)\xi-S(\xi,W)R(Y,Z)X &=& 0,
\end{eqnarray}
 for all vector fields $X,Y,Z,W\in TM$. Taking inner product of the above $(8.2)$ with the vector field $\xi$ and then replacing $W=\xi$ we obtain
 \begin{eqnarray}
  S(X,R(Y,Z)\xi)-S(\xi,R(Y,Z)\xi)\eta(X)+S(X,Y)\eta(R(\xi,Z)\xi) &&  \nonumber\\
  -S(\xi,Y)\eta(R(X,Z)\xi)+S(X,Z)\eta(R(Y,\xi)\xi)-S(\xi,Z)\eta(R(Y,X)\xi) &&  \nonumber\\
  +S(X,\xi)\eta(R(Y,Z)\xi)-S(\xi,\xi)\eta(R(Y,Z)X) &=& 0, \nonumber
\end{eqnarray}
 In view of $(2.11)$ and $(2.13)$ the above equation becomes
 \begin{eqnarray}
  (\alpha^2-\beta^2)[S(X,Y)\eta(Z)-S(X,Z)\eta(Y)-S(\xi,Y)\eta(X)\eta(Z) &&  \nonumber\\
  +S(\xi,Z)\eta(X)\eta(Y)]-S(\xi,\xi)\eta(R(Y,Z)X)&=& 0.
\end{eqnarray}
 Putting $Y=\xi$ in $(8.3)$ and then recalling $(2.12)$ we get
 \begin{eqnarray}
  (\alpha^2-\beta^2)[S(X,\xi)\eta(Z)-S(X,Z)-S(\xi,\xi)\eta(X)\eta(Z) &&  \nonumber\\
  +S(\xi,Z)\eta(X)-S(\xi,\xi)[g(X,Z)-\eta(X)\eta(Z)]]&=& 0.
\end{eqnarray}
 Using equations $(3.3)$ and $(3.4)$ in the above $(8.4)$ yields
 \begin{equation}
 (\alpha^2-\beta^2)[(r-2\lambda-2\mu+\beta)\eta(X)\eta(Z)-(r-2\lambda-\mu-\beta)g(X,Z)]=0.\nonumber
 \end{equation}
  Replacing $X=\xi$ in above we arrive at
 \begin{equation}
 (\alpha^2-\beta^2)(2\beta-\mu)\eta(X)=0,~~~~~\forall X\in TM.
 \end{equation}
Since for non-trivial case $\alpha^2\neq\beta^2$, from the above equation $(8.5)$ it follows that $\mu=2\beta$. Therefore in view of this and recalling $(3.6)$ we finally obtain $r=6\lambda+8\beta$. Therefore we can state the following
\begin{theorem}
  Let $(M,g)$ be a trans-Sasakian $3$-manifold admitting an $\eta$-Einstein soliton $(g,\xi,\lambda,\mu)$. If the manifold satisfies the curvature condition $S(\xi,X)\cdot R=0$, then it becomes a  manifold of constant scalar curvature $r=6\lambda+8\beta$.
\end{theorem}

\section{$\eta$-Einstein solitons on trans-Sasakian $3$-manifolds with torse-forming vector field}
This section is devoted to study the nature of $\eta$-Einstein solitons on trans-Sasakian $3$-manifolds with torse-forming vector field.
\begin{definition}
A vector field $V$ on a trans-Sasakian $3$-manifold is said to be torse-forming vector field \cite{Y1} if
\begin{equation}
\nabla_XV=fX+\gamma(X)V,
\end{equation}
where $f$ is a smooth function and $\gamma$ is a $1$-form.
\end{definition}
Now let $(g,\xi,\lambda,\mu)$ be an $\eta$-Einstein soliton on a trans-Sasakian $3$-manifold $(M,g)$ and assume that the Reeb vector field $\xi$ of the manifold is a torse-forming vector field. Then $\xi$ being a torse-forming vector field, by definiton from equation $(9.1)$ we have
\begin{equation}
\nabla_X\xi=fX+\gamma(X)\xi,
\end{equation}
$\forall X\in TM$, $f$ being a smooth function and $\gamma$ is a $1$-form.\\
Recalling the equation $(2.8)$ and taking inner product on both sides with $\xi$ we can write
\begin{equation}
g(\nabla_X\xi,\xi)=(\beta-1)\eta(X).
\end{equation}
Again from the equation $(9.2)$, applying inner product with $\xi$ we obtain
\begin{equation}
g(\nabla_X\xi,\xi)=f\eta(X)+\gamma(X).
\end{equation}
Combining $(9.3)$ and $(9.4)$ we get, $\gamma=(\beta-1-f)\eta$. Thus from $(9.2)$ it implies that, for torse-forming vector field $\xi$ in a trans-Sasakian $3$-manifold, we have
\begin{equation}
\nabla_X\xi=f(X-\eta(X)\xi)+(\beta-1)\eta(X)\xi.
\end{equation}
Now from the formula of Lie differentiation and using $(9.5)$ yields
\begin{eqnarray}
  (\mathcal{L}_{\xi}g)(X,Y) &=& g(\nabla_X\xi,Y)+g(\nabla_Y\xi,X) \nonumber\\
   &=& 2f[g(X,Y)-\eta(X)\eta(Y)]+2(\beta-1)\eta(X)\eta(Y).
\end{eqnarray}
Since $(g,\xi,\lambda,\mu)$ is an $\eta$-Einstein soliton, the equation $(1.2)$ holds. So in view of $(9.6)$, the equation $(1.2)$ reduces to
\begin{equation}
S(X,Y)=(\frac{r}{2}-\lambda+f)g(X,Y)+(f-\mu-\beta+1)\eta(X)\eta(Y).
\end{equation}
This implies that the manifold is an $\eta$-Einstein manifold. Again putting $Y=\xi$ in $(9.7)$ we get
\begin{equation}
S(X,\xi)=(\frac{r}{2}-\lambda-\mu-\beta+1)\eta(X).
\end{equation}
Combining $(9.8)$ with the equation $(2.15)$ implies
\begin{equation}
(\frac{r}{2}-\lambda-\mu-\beta+1)=2(\alpha^2-\beta^2).
\end{equation}
Again tracing out the equation $(9.7)$ we obtain
\begin{equation}
r=6\lambda+2\mu+4f+2\beta-2.
\end{equation}
Using the above equation $(9.10)$ in $(9.9)$, finally we get $\lambda=f-(\alpha^2-\beta^2)$. Therefore we have the following
\begin{theorem}
  Let $(g,\xi,\lambda,\mu)$ be an $\eta$-Einstein soliton on a trans-Sasakian $3$-manifold $(M,g)$, with torse-forming vector field $\xi$, then the manifold becomes an $\eta$-Einstein manifold and the soliton is expanding, steady or shrinking according as $f>(\alpha^2-\beta^2)$, $f=(\alpha^2-\beta^2)$ or, $f<(\alpha^2-\beta^2)$ respectively.
\end{theorem}

\medskip
\medskip
{\bf Acknowledgements.}
     The first author D. Ganguly is thankful to the National Board for Higher Mathematics (NBHM), India, for their financial support(Ref No: 0203/11/2017/RD-II/10440) to carry on this research work.

\bebi
\bibitem{Blg} A. M. Blaga,~{\em On gradient $\eta$-Einstein solitons},
~Kragujevac Journal of Mathematics,~vol-42(2),~pp.229-237,~(2018).
\bibitem{Bla} D. E. Blair,~{\em On the geometric meaning of the Bochner tensor},
~Geometriae Dedicata,~vol-4(1),~pp.33-38,~(1975).
\bibitem{Bla1} D. E. Blair,~{\em Riemannian geometry of contact and symplectic manifolds},
~Birkhauser,~Second Edition,~2010.
\bibitem{Boc} S. Bochner,~{\em Curvature and Betti numbers. $II$},
~Annals of Mathematics,~vol-50(2),~pp.77-93,~(1949).
\bibitem{Cr} C. Calin, M. Crasmareanu, ~{\em $\eta$-Ricci solitons on Hopf hypersurfaces in complex space forms},
~Revue Roumaine de Math. Pures et Appl.,~57(1),~pp.53-63,~(2012).
\bibitem{CK} J. T. Cho, M. Kimura, ~{\em Ricci solitons and real hypersurfaces in a complex space forms},
~Tohoku Math. J.,~61,~pp.205-212,~(2009).
\bibitem{CM} G. Catino, L. Mazzieri,~{\em Gradient Einstein solitons},
~Nonlinear Anal.,~132, ~pp.66-94,~(2016).
\bibitem{G}A. Gray,~{\em Einstein like manifolds which are not Einstein},
~Goem. Dedicata,~vol-7,~pp.259-280,~(1978).
\bibitem{GH}A. Gray, L. M. Hervella,~{\em The sixteen classes of almost Hermitian manifolds and their linear invariants},
~Ann. Math. Pura. Appl.,~vol-123(4),~pp.35-58,~(1980).
\bibitem{H1}R. S. Hamilton, ~{\em Three manifolds with positive Ricci curvature},
~Journal of Differential Geometry 17, ~255-306, ~(1982).
\bibitem{Ma} J. C. Marrero, ~{\em The local structure of trans-Sasakian manifolds},
~Ann. Math. Pura. Appl.,~vol-162(4),~pp.77-86,~(1992).
\bibitem{MC} M. Matsumoto, G. Ch\={u}man,~{\em On the C-Bochner curvature tensor},
~TRU mathematics,~vol-5,~pp.21-30,~(1969).
\bibitem{O}J. A. Oubina,~{\em New classes of almost contact metric structures},
~Publ. Math. Debrecen,~vol-32(4),~pp.187-193,~(1985).
\bibitem{PM}G. P. Pokhariyal, R. S. Mishra, ~{\em The curvature tensor and their relativistic significance},
    ~Yokohoma Math. J.,~vol-18,~pp-105-108,~(1970).
\bibitem{Sz}S. I. Szabo,~{\em Structure theorems on Riemannian spaces satisfying $R(X,Y)R=0$, $I$, The local version},~j. Differ. Geometry,~vol-17,~pp.531-582,~(1982).
\bibitem{YK}K. Yano, M. Kon,~{\em Structures on manifolds},~Series in Pure Mathematics,~vol.3,~(1984).
\bibitem{Y1}K. Yano,~{\em On torse-forming directions in Riemannian spaces},~Proceedings of the Imperial Academy Tokyo,~20,~pp.701-705,~(1944).
\eebi

\auas
    Dipen Ganguly\\
    Department~of~Mathematics, Jadavpur~University, Kolkata-700032, West Bengal,~India.\\
    E-mail: dipenganguly1@gmail.com\\

    \noindent Santu Dey\\
    Department~of~Mathematics, Bidhan~Chandra~College, Asansol-713304, West Bengal,~India.\\
    E-mail: santu.mathju@gmail.com\\

    \noindent Arindam Bhattacharyya\\
    Department~of~Mathematics, Jadavpur~University, Kolkata-700032, West Bengal,~India.\\
    E-mail: bhattachar1968@yahoo.co.in\\

\end{document}